\documentclass{article}
\usepackage{latexsym}
\usepackage[dvips]{graphicx}
\usepackage{amssymb}

\usepackage{amsmath}

\newcommand{\beqn}{\begin{equation}}
\newcommand{\eeqn}{\end{equation}}

\newcommand{\hur}{H_{\alpha}}
\newcommand{\whH}{\widehat{H}}

\newcommand{\pf}{\noindent{\sc Proof. }}
\newcommand{\epf}{\hfill$\Box$}
\newcommand{\mC}{\mathcal{C}}
\newcommand{\al}{\alpha}
\newcommand{\la}{\lambda}
\newcommand{\ld}{\ldots}
\newcommand{\pr}{\prime}
\newcommand{\f}{\frac}
\newcommand{\tf}{\tfrac}
\newcommand{\pa}{\partial}

\newtheorem{theorem}{Theorem}

\begin{document}

\title{A simple recurrence for covers of the sphere with\\
branch points of arbitrary ramification}
\author {I. P. Goulden\thanks{Supported by a Discovery Grant from
 NSERC. E-mail ipgoulden@uwaterloo.ca} ~ and Luis G. Serrano
\thanks{Research supported by a Postgraduate Scholarship from
NSERC. E-mail
 lgserrano@uwaterloo.ca}\\
University of Waterloo\\
Department of Combinatorics \& Optimization\\
Waterloo, Ontario N2L 3G1, Canada}
\maketitle

\begin{abstract}
The problem of counting ramified covers of a Riemann surface up to homeomorphism
was proposed by Hurwitz in the late 1800's. This problem translates combinatorially
into factoring a permutation of specified cycle type, with certain conditions
on the cycle types of the factors, such as minimality and transitivity.

Goulden and Jackson have given a proof for the number of minimal, transitive
factorizations of a permutation into transpositions. This proof involves a
partial differential equation for the generating series, called the Join-Cut
equation.
Recently, Bousquet-M\'{e}lou and Schaeffer have found the number of minimal, transitive
factorizations of a permutation into arbitrary unspecified factors. This was proved by a
purely combinatorial argument, based on a direct bijection between factorizations
and certain objects called $m$-Eulerian trees.

In this paper, we give a simple partial differential equation for
Bousquet-M\'{e}lou and Schaeffer's generating series, and for
Goulden and Jackson's generating series, as well as a new proof of
the result by Bousquet-M\'{e}lou and Schaeffer. We apply algebraic
methods based on Lagrange's theorem, and combinatorial methods based
on a new use of Bousquet-M\'{e}lou and Schaeffer's $m$-Eulerian trees.
\end{abstract}

\section{Introduction}\label{intro}

\subsection{Background and notation}

The subject of this paper is a mathematical problem that has attracted significant
attention in the last decade, with contributions from researchers in algebraic
combinatorics, algebraic geometry, and mathematical physics. Its origins are in
the work of Hurwitz~\cite{H}, dating from the 1890's.

First, for some notation.
A {\em partition} is a weakly ordered list of positive integers $\la =(\la_1,\ld ,\la_k)$,
where $\la_1\ge\ld\ge\la_k$.
The integers $\la_1,\ld ,\la_k$ are called the {\em parts} of the partition $\la$,
and we denote the number of parts by $l(\la )=k$. If $\la_1+\ld +\la_k=n$, then $\la$ is
a partition of $n$, and we write $\la\vdash n$. In the symmetric group $S_n$ on $\{ 1,\ld ,n\}$,
let $\mC_{\la}$ be the conjugacy class consisting of all permutations whose disjoint cycle lengths
are specified by the parts of the partition $\la$. If $\la$ has $d_i$ parts equal
to $i$, for $i\ge 1$, then it is well known
that $\vert\mC_{\la}\vert=n!/\prod_{i\ge 1}i^{d_i}d_i!$.

In this paper, we consider branched covers of the sphere by an $n$-sheeted
Riemann surface of genus~$g$. Suppose that the branch points
are $P_0, P_1,\ld ,P_k$, with ramification at $P_i$ specified by $\pi_i\in P_i$,
for $i=0,\ld ,k$. (This means that if one walks in a small neighbourhood, counterclockwise,
around $P_i$, starting at sheet $j$, then one ends at sheet $\pi_i(j)$.)
The {\em monodromy} condition states that, for consistency, we must have
\begin{equation}\label{factn}
\pi_1\ld\pi_k=\pi,
\end{equation}
where $\pi=\pi_0^{-1}$. For the cover to be connected, one must be able to move
from sheet $j$ to sheet $m$, by walking around suitable branch points, for
any choices of $j$ and $m$. This means that the group generated by $\pi_1,\ld ,\pi_k$ must
be a transitive subgroup of $S_n$. In this case, we call the factorization
of $\pi$ in~(\ref{factn}) a {\em transitive} factorization. Finally,
if $\pi_i\in\mC_{\al_i}$, $i=0,\ld ,n$, then we say that the {\em ramification-type}
of $P_i$ is $\al_i$, and
the {\em Riemann-Hurwitz} formula then gives
\begin{equation}\label{RHf}
\sum_{i=0}^k (n-l(\al_i ))=2n-2+2g.
\end{equation}
In the case $g=0$, we call the factorization a {\em minimal} transitive factorization.
Hurwitz~\cite{H} proved that the ordered $\pi_i$'s specify a unique cover
up to homeomorphism, so we can count (homeomorphically) distinct branched
covers with specified ramification by counting transitive factorizations.
A branch point whose ramification is a transposition is called
a {\em simple} branch point.

In genus $0$, there are two explicit, general results known. For the first result,
let $H_{\al}$ be the number of distinct covers of the sphere by an $n$-sheeted Riemann
surface of genus $0$, with fixed ramification of type $\al$ at $P_0$, and
all other branch points are simple. Then Goulden and Jackson~\cite{GJ}
(see also Hurwitz~\cite{H2}) have proved that
\beqn\label{HN}
 \hur= n^{l(\alpha)-3} (n + l(\alpha) - 2)! \prod_{i \ge 1}
 \left(\frac{i^i}{(i-1)!}\right)^{d_i},
\eeqn
where $\al$ has $d_i$ parts equal to $i$, $i\ge 1$. In this case, the
Riemann-Hurwitz formula~(\ref{RHf}) gives the number of simple
branch points as $n+l(\al )-2$.

For the second explicit result, let $G_{\al}(m)$ be the number of distinct covers of the
sphere by an $n$-sheeted Riemann
surface of genus $0$, with fixed ramification of type $\al$ at $P_0$, and
$m$ other branch points with arbitrary ramification. Then
Bousquet-M\'{e}lou and Schaeffer~\cite{BMS} have proved that
\beqn\label{GN}
G_{\alpha}(m) = m \frac{\{ (m-1)n-1\} !}{ \{ (m-1)n - l(\alpha) + 2 \} !}
\prod_{i \ge 1}\left\{ i{mi-1 \choose i}\right\}^{d_i},
\eeqn
where $\al$ has $d_i$ parts equal to $i$, $i\ge 1$. In this case, the
Riemann-Hurwitz formula~(\ref{RHf}) gives the total number
of cycles in the arbitrary factors as $c(\al )=(m-1)n-l(\al )+2$.

The following version of Lagrange's Theorem, as it appears
in~\cite{GJBook}, is used in this paper.
We write $[A]B$ to mean ``the coefficient of $A$ in $B$''.

\begin{theorem}\label{lagrange}\textbf{(Lagrange)}
Let $\phi(\lambda) \in R[[\lambda]]_1$. Then there exists a unique
formal power series $w(t) \in R[[t]]_0$  such that $w = t \phi(w)$. Moreover,
if $f(\lambda) \in R((\lambda))$, then for $n \ne 0$,
$$[t^n]f(w) = \frac{1}{n}[\lambda^{n-1}] f^{'}(\lambda) \phi^n (\lambda).$$
\end{theorem}

\subsection{Outline}

In Section~\ref{pdesection}, we apply Lagrange's Theorem to prove that
a generating series for Bousquet-M\'{e}lou and Schaeffer's~\cite{BMS}
numbers $G_{\al}(m)$, given in~(\ref{GN}), satisfies a simple
quadratic partial differential equation. This equation
is given in Theorem~\ref{Hquad}. In Section~\ref{Etrees}, we
give a combinatorial proof of Theorem~\ref{Hquad}, and the
associated Lagrangian equations, using a minor adaptation of
the combinatorial objects introduced by
Bousquet-M\'{e}lou and Schaeffer~\cite{BMS} in their proof
of~(\ref{GN}). In Section~\ref{Simpbranch}, we state, as Theorem~\ref{kdv2},
a simple quadratric partial differential equation for a generating
series for the numbers $H_{\al}$, given in~(\ref{HN}).
We know of no direct combinatorial or geometric proof of either
Theorem~\ref{Hquad} or Theorem~\ref{kdv2}, despite the
apparent simplicity of these equations.

Moreover, the numbers $H_{\al}$ satisfy a well-known partial differential
equation called the Join-Cut equation. (This equation was
used in~\cite{GJ} to obtain~(\ref{HN}).) The Join-Cut equation
has a simple combinatorial proof, by analyzing the transpositions
in a corresponding factorization, and a simple geometric proof,
by analyzing branch points. The Join-Cut equation
also extends in a straightforward way to the case of
arbitrary genus. In arbitrary genus, the corresponding numbers
(in which one fixed branch point has arbitrary ramification,
and all others are simple) are given by the
ELSV formula~\cite{ELSV}, as a Hodge integral. However, we
know of no simple combinatorial or geometric proof of
Theorem~\ref{kdv2}.

For the numbers $G_{\al}(m)$, we know of no associated results analogous
to those described above for the numbers $H_{\al}$.
We speculate that analogues of some of these  associated results
exist, and are particularly interested in extensions to higher genus.

\section{Branch points of arbitrary ramification}\label{pdesection}

In this section we consider the following generating series for
the numbers $G_{\al}(m)$ given in~(\ref{GN}):
\begin{equation}\label{Ggs}
G \equiv G(z,u,x, p_1, p_2, \ld)
= \sum_{n \ge 1} \sum_{\alpha \vdash n} G_{\alpha}(m) \frac{z^n}{n!}
 \vert{\mC}_{\al}\vert u^{l(\alpha)} p_{\al} x^{(m-1)n - l(\al) + 2}.
\end{equation}
Here $z$ is an exponential indeterminate
marking the elements in $S_n$, $u$ is an ordinary indeterminate marking the cycles
in $\pi$, and $x$ is an ordinary indeterminate marking the total
number of cycles of the factors $\pi_1, \ldots, \pi_m$, which by the
minimality condition is $c(\alpha)=(m-1)n - l(\alpha) + 2$.

For fixed $m\ge 1$,
let $A(t) = \sum_{i \ge 1} {mi-1 \choose i} u p_i t^i$,
and $w\equiv w(z,u,x, p_1, p_2, \ld)$ be the
unique formal power series satisfying
\begin{equation}\label{lagw}
w = z (x+A(w))^{m-1}.
\end{equation}
In the following result, we use~(\ref{GN}) to express two partial
derivatives of $G(z,u,p_1,p_2,\ld )$ in terms of $A(w)$.
Of course, we can determine $w$  and $A(w)$ by Lagrange's Theorem.

\begin{theorem}\label{mainequation}
For $G$ given by~(\ref{Ggs}), and $w$ given by~(\ref{lagw}), we have
\begin{equation}\label{main1}
\left(z \frac{\partial}{\partial z} \right) \left((m-1)
z \frac{\partial}{\partial z} + 1\right) G = \frac{m}{m-1} x\, A(w) + \frac{m}{2(m-1)} A(w)^2,
\end{equation}
and
\begin{equation}\label{main2}
%z \frac{\pa}{\pa z} \left((m-1)
%z \frac{\pa}{\pa z} - u\f{\pa}{\pa u}+ 2\right) G = \frac{m}{m-1}A(w).
\left(z \f{\pa}{\pa z} \right) \left( \f{\pa}{\pa x} \right) G = \frac{m}{m-1}A(w).
\end{equation}
\end{theorem}

\pf
From~(\ref{GN}), and the fact that $\vert{\mC}_{\al}\vert=n!/{\prod_{i \ge 1}i^{d_i}d_i!}$, we
obtain
\begin{equation*}
G=\sum_{n \ge 1}z^n \sum_{\alpha \vdash n}m \f{\{ (m-1)n-1\} !}{ \{ (m-1)n - l(\alpha) + 2 \} !}
x^{(m-1)n - l(\al) + 2} \prod_{i \ge 1} {mi-1 \choose i}^{d_i}\f{p_i^{d_i}}{d_i !} u^{d_i},
\end{equation*}
where $\al$ has $d_i$ parts equal to $i$, for $i\geq 1$. Thus, we have, for $n\geq 1$,
\begin{eqnarray*}
[z^n]G & = & \sum_{k \ge 1}m \f{\{ (m-1)n-1\} !}{ \{ (m-1)n -k+ 2 \} !}
 \sum_{\stackrel{d_1 + 2d_2 + \cdots = n}{d_1 + d_2 + \cdots = k}}
x^{(m-1)n - k + 2}
\prod_{i \ge 1} {mi-1 \choose i}^{d_i}\f{p_i^{d_i}}{d_i !} u^{d_i} \\
&=& \sum_{k \ge 1}m \f{\{ (m-1)n-1\} !}{ \{ (m-1)n -k+ 2 \} !}
[\la^n]\f{A(\la )^k}{k!} x^{(m-1)n - k + 2} ,
\end{eqnarray*}
from the multinomial theorem. Thus, we obtain
\begin{eqnarray*}
[z^n]z \f{\pa}{\pa z} \left((m-1)z \f{\pa}{\pa z} + 1\right) G
&=&\sum_{k \ge 1} n\{(m-1)n+ 1\}m \f{\{ (m-1)n-1\} !}{ \{ (m-1)n -k+ 2 \} !}
[\la^n]\f{A(\la )^k}{k!} x^{(m-1)n - k + 2} \\
&=& \sum_{k \ge 1}\f{m}{(m-1)\{ (m-1)n+2\} }{(m-1)n+2 \choose k} [\la^n] A(\la )^k x^{(m-1)n - k + 2} \\
&=&\f{m}{(m-1)\{ (m-1)n+2\} } [\la^n](x+A(\la ))^{(m-1)n+2}\\
&=&\f{m}{(m-1)\{ (m-1)n+2\} }\f{1}{n}[\la^{n-1}]\f{d}{d \la}(x+A(\la ))^{(m-1)n+2}\\
&=&\f{m}{m-1}\f{1}{n}[\la^{n-1}] A^{\pr}(\la)(x+A(\la ))^{(m-1)n+1}\\
&=&\f{m}{m-1}[z^n]\left( x A(w)+\tf{1}{2}A(w)^2\right),
\end{eqnarray*}
from Theorem~\ref{lagrange}, and~(\ref{main1}) follows, since both $G$ and $A(w)$ are
equal to $0$ at $z=0$.

Similarly, we obtain
\begin{eqnarray}
[z^n]z \frac{\pa}{\pa z} \left(
\f{\pa}{\pa x} \right) G
&=&\sum_{k \ge 1} n\{(m-1)n-k+2\}m \f{\{ (m-1)n-1\} !}{ \{ (m-1)n -k+ 2 \} !}
[\la^n]\f{A(\la )^k}{k!} x^{(m-1)n - k + 1} \nonumber \\
&=& \sum_{k \ge 1}\f{m}{(m-1)\{ (m-1)n+1\} }{(m-1)n+1 \choose k}[\la^n]A(\la )^ k x^{(m-1)n - k + 1} \nonumber \\
&=&\f{m}{(m-1)\{ (m-1)n+1\} } [\la^n](x+A(\la ))^{(m-1)n+1} \nonumber \\
&=&\f{m}{(m-1)\{ (m-1)n+1\} }\f{1}{n}[\la^{n-1}]\f{d}{d \la}(x+A(\la ))^{(m-1)n+1} \nonumber \\
&=&\f{m}{m-1}\f{1}{n}[\la^{n-1}] A^{\pr}(\la)(x+A(\la ))^{(m-1)n} \nonumber \\
&=&\f{m}{m-1}[z^n]A(w), \label{GAw}
\end{eqnarray}
from Theorem~\ref{lagrange}, and~(\ref{main2}) follows.
\epf
\bigskip

Note that since
$$\left( x \f{\pa}{\pa x} \right) G = \left( (m-1) z \f{\pa}{\pa z} - u \f{\pa}{\pa u} + 2 \right) G,$$
then we can eliminate $A(w)$ between equations~(\ref{main1}) and~(\ref{main2}), to
obtain a quadratic partial differential equation for $G$, given in the following result.

\begin{theorem}\label{Hquad}
The generating series $G$ satisfies the partial differential equation
\begin{equation}\label{recG}
2m z\f{\pa}{\pa z} \left(u\f{\pa}{\pa u} - 1\right) G =
(m-1)\left(z\f{\pa}{\pa z} \f{\pa}{\pa x} G\right)^2.
\end{equation}
\end{theorem}

\section{$m$-Eulerian trees}\label{Etrees}

In this section we will provide a combinatorial interpretation for the partial
differential equations in Section \ref{pdesection}.
With the purpose of interpreting~(\ref{recG}), we will make use of certain objects
called $m$-Eulerian trees, defined by Bousquet-M\'{e}lou and Schaeffer in~\cite{BMS}.
A {\em leaf} in a tree is a vertex of degree $1$, and an {\em inner vertex}
is one that is not a leaf. The {\em inner degree} of a vertex
is the number of inner neighbours.
An $m$-{\em Eulerian} tree is a (properly) coloured tree where the colours of the
vertices are black and white, and has the following two properties:
\begin{itemize}
\item Every inner black vertex has total degree $m$ and inner
degree $1$ or $2$,
\item Every inner white vertex has total degree $mi$ for
some $i \ge 1$, and exactly $i-1$ neighbours of inner degree $1$.
\end{itemize}

Now, it is straightforward to prove by a counting argument that the number of
black leaves is $m$ greater than the number of white leaves in an $m$-Eulerian tree.
Bousquet-M\'{e}lou and Schaeffer~\cite{BMS} have proved that there is a
unique way of inserting edges between white leaves 
and black leaves in such a  way that the outer face in the resulting planar map
contains the $m$ unpaired black leaves.
An $m$-Eulerian tree is said to be {\em planted} if a black leaf is specified as
a root vertex.
Furthermore, it is said to be \em balanced \em if the root is precisely one of
these $m$ black unpaired leaves. Note that in this paper, we allow
an $m$-Eulerian
tree that is not planted (notice
the slight difference between this definition and that in \cite{BMS}).

\begin{figure}[ht]
\centerline{ \includegraphics[scale= .50]{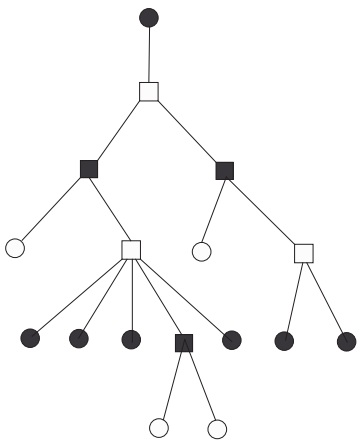} }
\caption{A Balanced m-Eulerian Tree\label{euleriantree} }
\end{figure}

If $\alpha \vdash n$ has $d_i$ parts of size $i, i \ge 1$, an $m$-Eulerian tree
of {\em type} $\alpha$ is defined as
having $n-1$ black inner vertices, $l(\alpha)$ white
inner vertices, $d_i$ of which have total degree $mi$, and $c(\alpha)$ black
leaves. An example can be seen in Figure \ref{euleriantree}, where
inner vertices  are drawn
as squares, and leaves are drawn as circles.

Bousquet-M\'{e}lou and Schaeffer have proved~(\ref{GN}) by finding a bijection
between minimal transitive factorizations of a permutation of
type $\alpha$ into $m$ factors, and balanced $m$-Eulerian trees of
type $\alpha$.
In this bijection, the element $1$ becomes a new separate black vertex of degree $m$, that is 
%In this bijection, a particular element from $\{1, \ld ,n\}$ is
%selected, and becomes the label attached to a new vertex of degree $m$, that is
inserted into the outer face of the map described above, and joined
to the $m$ unpaired black leaves.
The remaining elements
of $\{2, \ld ,n\}$ correspond to the inner black vertices, the cycles of $\alpha$ to
the inner white vertices, and the cycles of the $m$ factors correspond to the
black leaves.  Furthermore, to turn the $m$-Eulerian tree into a balanced tree,
one can plant it in $m$ ways. Therefore, if $T_{\alpha}(m)$ is the number
of $m$-Eulerian trees of type $\alpha$, then we have
$$mT_{\alpha}(m) = \f{\vert \mathcal{C}_{\alpha} \vert}{(n-1)!} G_{\alpha}(m).$$
Thus, by defining the generating series $T$ as
$$T \equiv T(z,u,x, p_1, p_2, \ldots)
= \sum_{n \ge 1} \sum_{\alpha \vdash n} T_{\alpha}(m) z^{n-1} u^{l(\alpha)} p_{\alpha} x^{(m-1)n - l(\alpha) + 2},$$
where $z$ is an ordinary marker for the inner black vertices, $u$ an ordinary marker for the inner white vertices, and $x$ an ordinary marker for the black leaves, we obtain
\beqn\label{TtoG}
mT = \f{\pa}{\pa z} G.
\eeqn
Applying~(\ref{TtoG}), equation~(\ref{recG}) for $G$ becomes the following
equation for T:
\beqn\label{recT}
2 \left( u \f{\pa}{\pa u} - 1 \right) T = 
(m-1) z \left( \f{\pa}{\pa x} T \right)^2.
\eeqn
Using $m$-Eulerian trees, we have the following combinatorial proof of~(\ref{recT}).
Consider two $m$-Eulerian trees (possibly identical). Root each one at a black
leaf (not necessarily an unpaired one). Now, identify these two roots,
to become a new inner black vertex of degree $2$, thus obtaining the generating
series $z \left(\frac{\partial}{\partial x} T\right)^2$. Note that the tree obtained
has the properties of an $m$-Eulerian tree at every vertex except for the new one, and
all we need to do to satisfy the conditions of an $m$-Eulerian tree is to
attach $m-2$ leaves to this new vertex. This can be done in $m-1$ ways, since
we simply decide how many leaves go on each side of the vertex,
thus giving the right hand side of~(\ref{recT}).

\begin{figure}[ht]
\centerline{ \includegraphics[scale= .50]{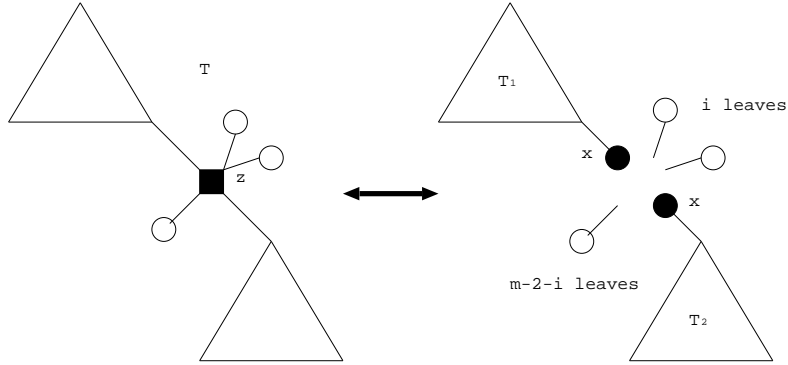} }
\caption{Cutting the Tree\label{cuttree} }
\end{figure}

For the left hand side of~(\ref{recT}), note that this process determines a
new $m$-Eulerian tree, with a selected inner black vertex of inner degree $2$.
In order to count the number of inner black vertices of inner degree $2$,
root the $m$-Eulerian tree at any black leaf. Then every inner white vertex,
except for the one adjacent to the root, can be paired up uniquely with the
inner black vertex of inner degree $2$ that immediately follows it
in the unique path from the white vertex to the root. Therefore
there are $l(\alpha) - 1$ of these black vertices. Thus 
the generating series for selecting one of these
is $\left( u\frac{\partial}{\partial u} -1\right)T$. Note that in this bijection
every Eulerian tree is constructed twice, which accounts for the factor
of $2$ on the left hand side of~(\ref{recT}). This process, which is clearly
bijective, is illustrated in Figure \ref{cuttree}.
\begin{figure}[ht]
\centerline{ \includegraphics[scale= .50]{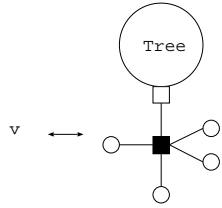} }
\caption{A pseudo-Eulerian Tree\label{pseudo} }
\end{figure}

We can also
interpret the Lagrangian equations in Section~\ref{pdesection}.
In order to do so, we define
a \em pseudo-Eulerian tree \em. This tree is constructed by removing
the root of a planted $m$-Eulerian tree (note that it does not need to be balanced), and adding an inner black vertex
in its place, with $m-1$ new white leaves emanating from it. As a final step,
we root this new tree at one of these $m-1$ white leaves.
A pseudo-Eulerian tree is illustrated in Figure~\ref{pseudo}.
For $\al \vdash n$,
let $v_{\alpha}$ be the number of
pseudo-Eulerian trees with $n$ inner black vertices, $l(\alpha)$ inner
white vertices, and $(m-1)n - l(\al) + 1$ black leaves, and
define $v$ as the generating series
$$v \equiv v(z, u, x, p_1, p_2, \ldots) =
 \sum_{\alpha \vdash n} v_{\alpha} z^n u^{l(\alpha)} p_{\alpha} x^{(m-1)n - l(\al) + 1}.$$

\clearpage

\noindent
Then, by the construction described above,
\begin{equation}\label{defv}
v = (m-1) z \f{\pa}{\pa x} T.
\end{equation}
Note that by (\ref{main2}), (\ref{TtoG}), and (\ref{defv}), we can conclude that 
\begin{equation}\label{vAw}
v = A(w).
\end{equation}

\begin{figure}[ht]
\centerline{ \includegraphics[scale= .50]{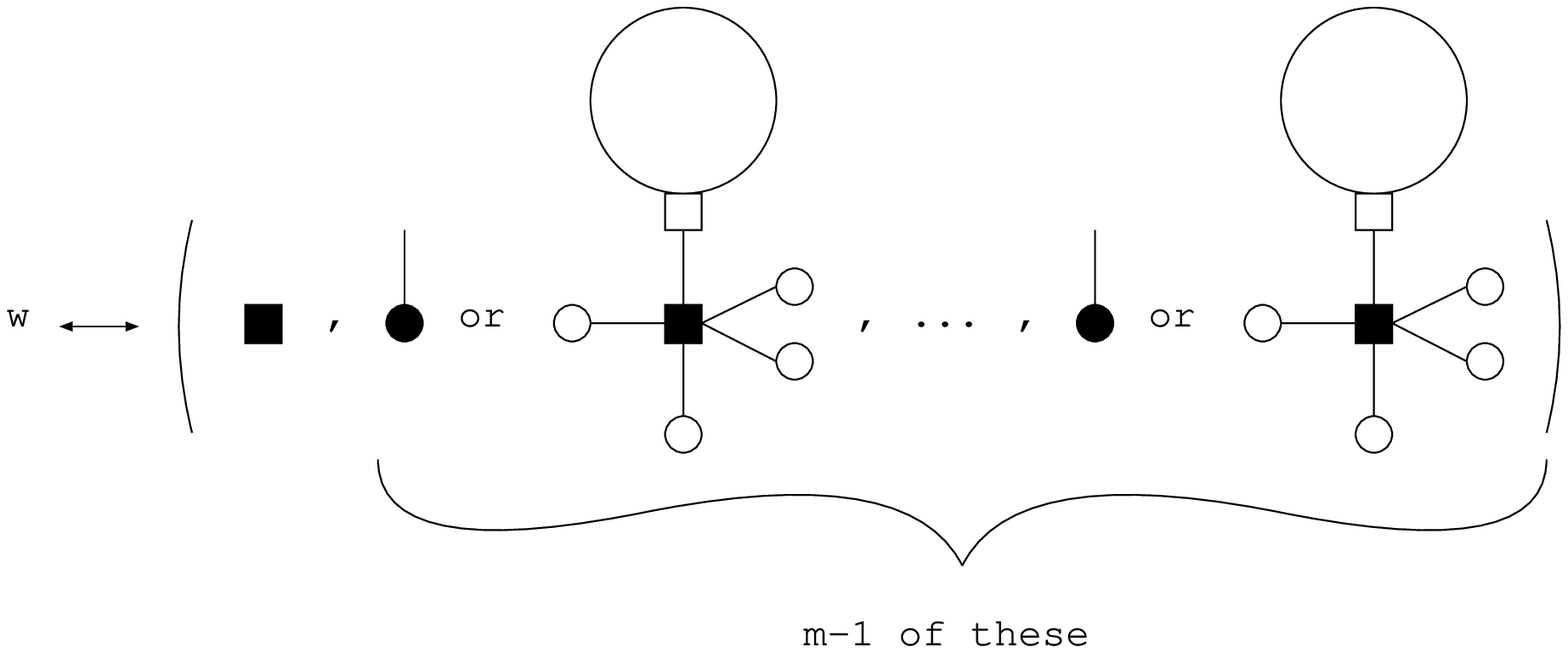} }
\caption{A list\label{list} }
\end{figure}

Now, define a \em list \em as a list of $m$ objects, the first of which is an
inner black vertex, and the remaining $m-1$ are a choice of either a black leaf,
or a pseudo-Eulerian tree, as seen on Figure \ref{list}. Note that this implies
that the generating series for lists, where the variables $z$, $x$, $u$,
and $p_i$ mark the same objects as before, is $z(x+v)^{m-1}$, which
by (\ref{lagw}) and (\ref{vAw}) is equal to $w$.

\begin{figure}[ht]
\centerline{ \includegraphics[scale= .50]{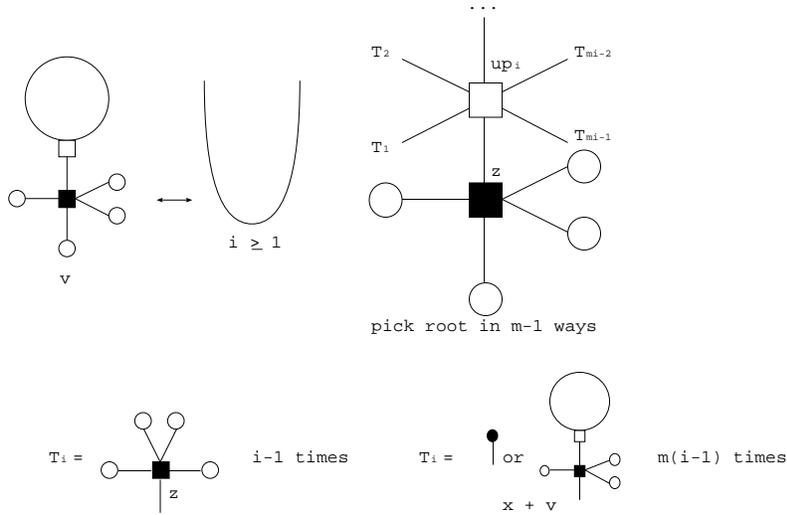} }
\caption{Decomposition of pseudo-Eulerian trees\label{pseudodecomposition} }
\end{figure}

In order to interpret (\ref{vAw}), we decompose a pseudo-Eulerian tree.
Note that a pseudo-Eulerian planted tree is planted at a white leaf, which is
selected in $m-1$ ways, and adjacent to a black inner vertex of inner
degree $1$, marked by $z$. This inner vertex is adjacent to exactly one inner white
vertex of degree $mi$, for some $i \ge 1$, marked by $up_i$
Now, for the remaining $mi-1$ edges that emanate from this white
vertex, $i-1$ of them must be adjacent to an inner black vertex of inner
degree $1$, and their positions can be chosen in ${mi-1 \choose i-1}$ ways. Each one of
them is marked by $z$. So we have an ordered set of $i$ copies of the
variable $z$, starting from the black vertex adjacent to the root, and going
clockwise around the white vertex specified above. Now, on the
remaining $mi-i$ edges adjacent to this white vertex, we must attach either
a black leaf, or a pseudo-Eulerian tree, together giving a factor of $x+v$. Again, taking them clockwise, we
have an ordered list of $(m-1)i$ objects, each one being a choice of a
black leaf or a pseudo-Eulerian tree. Note that all these, together with
the $m$ inner black vertices of inner degree one above, can be arranged
into an ordered set of $i$ lists, in the canonical way, each marked
by $w$. The equation obtained by this reasoning is $v = (m-1) \sum_{i \ge 1} {mi-1 \choose i-1} u p_i w^i = A(w)$. This process, which is clearly reversible, can be seen more
clearly in Figure \ref{pseudodecomposition}.

This concludes the combinatorial interpretation of the partial differential
equations in Section \ref{pdesection}. Moreover, by
Theorem~\ref{lagrange},~(\ref{lagw}), and~(\ref{defv}), one can recover~(\ref{GN}),
thus finding a different proof of the result in \cite{BMS}.

\section{Simple branch points}\label{Simpbranch}

In this section we consider the following generating series for
the numbers $H_{\al}$ given in~(\ref{HN}):
\begin{equation}\label{Hgs}
H \equiv H(z,p_1,p_2,\ld )
= \sum_{n \ge 1} \sum_{\alpha \vdash n} \hur\frac{z^n}{n!}
 \f{\vert{\mC}_{\al}\vert}{( n + l(\alpha) - 2) !} p_{\alpha}.
\end{equation}
Let $B(t)=\sum_{i\ge 1}\f{i^i}{i!} p_i t^i$, and let $s\equiv s(z,p_1,p_2,\ld )$ be the unique
formal power series satisfying
\begin{equation}\label{lagrangeeq}
s = z \exp B(s).
\end{equation}
The following result, from~\cite{GJ},  expresses three partial
derivatives of $H(z,p_1,p_2,\ld )$ in terms of $s$ and $B(s)$.
\begin{theorem}\label{bigprop}
For $H$ given by~(\ref{Hgs}), and $s$ given by~(\ref{lagrangeeq}), we have
\begin{equation}\label{prop1}
\left( z \frac{\partial}{\partial z}\right)^2 H = B(s),
\end{equation}
\begin{equation}\label{prop2}
z \f{\pa }{\pa z}H = \sum_{i \ge 1} \f{i^{i-1}}{i!} p_i s^i - \tf{1}{2}B(s)^2,
\end{equation}
and for $k \ge 1$,
\begin{equation}\label{prop3}
z \f{\pa}{\pa z}\f{\pa}{\pa p_k}H = \frac{k^{k-1}}{k!} s^k.
\end{equation}
\end{theorem}

The following result
gives a simple quadratic partial differential equation for $H$, that follows
immediately from Theorem~\ref{bigprop}. This result has been obtained
by Goulden and Jackson (private communication).

\begin{theorem}\label{kdv}
The generating series $H$ satisfies the partial differential equation
$$z \f{\pa }{\pa z}H = z\f{\pa}{\pa z}\sum_{k \ge 1}  p_k \f{\pa }{\pa p_k }H
-\tf{1}{2}\left(\left(z\f{\pa}{\pa z}\right)^2H\right)^2.$$
\end{theorem}

\pf Substitute (\ref{prop1}) and (\ref{prop3}), after adding over all $k$,
into~(\ref{prop2}).
\epf
\bigskip

This partial differential equation can be transformed in a straightforward way to have
a simpler form, by considering the series
$$\whH=
\sum_{n \ge 1} \sum_{\alpha \vdash n} \hur \frac{z^n}{n!}
\f{ \vert{\mC}_{\al}\vert}{(n-l(\al )+2)!} p_{\alpha} u^{l(\al)}.$$
Note that, in the generating series $\whH$, $z$ is an exponential indeterminate
marking the elements in $S_n$, u is an ordinary indeterminate marking the cycles
in $\pi$, and $p_i$ is an ordinary indeterminate marking the cycles of length $i$ in $\pi$.

\begin{theorem}\label{kdv2}
The generating series $\whH$ satisfies the partial differential equation
\begin{equation*}
z \f{\pa }{\pa z}\whH = z\f{\pa}{\pa z}u\f{\pa}{\pa u}\whH
-\tf{1}{2}\left(\left(z\f{\pa}{\pa z}\right)^2\whH\right)^2.
\end{equation*}
\end{theorem}

\pf
The result follows from Theorem~\ref{kdv} and the fact
that $\left( u\f{\pa}{\pa u} \right) \whH =
 \left( \sum_{k \ge 1} p_k \f{\pa}{\pa p_k} \right) \whH.$  %~(\ref{pdlem}).
\epf
\bigskip

Note that in this context, $u$ keeps track of the cycles in $\pi$. The apparent simplicity
of Theorem~\ref{kdv2}, especially as it involves only the two variables $z$ and $u$,
suggests that there may be a simple combinatorial interpretation.
We have looked for such a combinatorial interpretation, but have been unable to find one.


\begin{thebibliography}{9}


       \bibitem{BMS}
       M. Bousquet-M\'{e}lou, G. Schaeffer,
       \emph{Enumeration of Planar Constellations},
       Adv. Appl. Math. \textbf{24} (2000), 337--368.


%      \bibitem{C}
%      Cedric Chauve,
%      \emph{A Bijection Between Planar Constellations and some Colored Lagrangian Trees},
%      Discrete Mathematics and Theoretical Computer Science \textbf{6} (2003), 13--40.


       \bibitem{ELSV}
       T. Ekedahl, S. Lando, M. Shapiro, A. Vainshtein,
       \emph{Hurwitz Numbers and Intersections on Moduli Spaces of Curves},
       Invent. Math. \textbf{146} (2001), 297--327.


%      \bibitem{F}
%      O. Forster,
%      \emph{``Lectures on Riemann Surfaces''},
%      Springer-Verlag, New York, 1981.


       \bibitem{GJBook}
       I.P. Goulden, D.M. Jackson,
       \emph{``Combinatorial Enumeration''},
       John Wiley \& Sons, New York, 1983 (Dover Reprint, 2004).


       \bibitem{GJ}
       I.P. Goulden, D.M. Jackson,
       \emph{Transitive Factorizations into Transpositions and Holomorphic Mappings on the Sphere},
       Proc. Amer. Math. Soc. \textbf{125} (1997), 51--60.


       \bibitem{H}
       A. Hurwitz,
       \emph{\H{U}ber Riemann'sche Fl\H{a}chen mit gegebenen Verzweignungspunkten},
       Math. Ann. \textbf{39} (1891), 1--60.

	\bibitem{H2}
       A. Hurwitz,
       \emph{\H{U}ber die Anzahl der Riemann'schen Fl\H{a}chen mit
 gegebenen Verzweignungspunkten},
       Math. Ann. \textbf{55} (1902), 53--66.

\end{thebibliography}
\end{document}